\documentclass[a4paper,12pt,notitlepage]{amsart}

\usepackage[utf8]{inputenc}
\usepackage[spanish,english]{babel}

\usepackage[T1]{fontenc} 
\usepackage{lmodern} 

\usepackage[
    margin=28mm,
    marginparwidth=25mm,
    marginparsep=2mm,
    bottom=25mm
    ]{geometry}

\usepackage{amsmath}
\usepackage{amsfonts}
\usepackage{amssymb}
\usepackage{amsthm}

\usepackage{xcolor}

\usepackage{tcolorbox}

\usepackage{caption}

\usepackage[pdfencoding=unicode, psdextra]{hyperref}

\usepackage{cleveref}

\usepackage[shortlabels]{enumitem}

\usepackage{mathtools}

\usepackage{csquotes}

\usepackage[sortcites,backend=biber, style=numeric, maxbibnames=10,citestyle=numeric-comp]{biblatex} 
\DeclareNameAlias{author}{family-given} 
\DeclareFieldFormat[misc]{title}{``#1''} 
\newtheorem{theorem}{Theorem}[section]
\newtheorem{corollary}[theorem]{Corollary}
\newtheorem{proposition}[theorem]{Proposition}
\newtheorem{lemma}[theorem]{Lemma}
\newtheorem{remark}[theorem]{Remark}
\newtheorem{definition}[theorem]{Definition}

\newcommand\abs[1]{\left|#1\right|}
\newcommand\norm[1]{\left\Vert#1\right\Vert}
\newcommand\floor[1]{\left\lfloor#1\right\rfloor}

\newcommand{\NN}{\mathbb N}

\newcommand{\RR}{\mathbb R}
\newcommand{\CC}{\mathbb C}

\newcommand{\BB}{\mathcal B}
\newcommand{\DD}{\mathbb D}
\newcommand{\TT}{\mathbb T}
\newcommand{\al}{\alpha}
\newcommand{\be}{\beta}
\newcommand{\e}{\varepsilon}

\DeclareMathOperator*{\esssup}{ess\,sup}

\definecolor{newred}{HTML}{ff6961}
\definecolor{neworange}{HTML}{ffb480}
\definecolor{newyellow}{HTML}{f8f38d}
\definecolor{newgreen}{HTML}{42d6a4}
\definecolor{newcyan}{HTML}{08cad1}
\definecolor{newblue}{HTML}{59adf6}
\definecolor{newviolet}{HTML}{9d94ff}
\definecolor{newpurple}{HTML}{c780e8}

\numberwithin{equation}{section}

\begin{document}

\title[A Ritt-Kreiss condition: Spectral localization and norm estimates]
{A Ritt-Kreiss condition: Spectral localization and norm estimates}

\author[A. Mahillo]{Alejandro Mahillo}
\address{\newline
    Alejandro Mahillo \newline
    Departamento de Matemáticas,
    Instituto Universitario de Matemáticas y Aplicaciones,
    Universidad de Zaragoza,
    50009 Zaragoza, Spain.}
\email{almahill@unizar.es}

\author[S. Rueda]{Silvia Rueda}
\address{\newline
    Silvia Rueda \newline
   Departamento de Matemática, Facultad de Ciencias,
    Universidad del Bío Bío,
    Concepción, Chile.}
\email{sruedas@ubiobio.cl}

\thanks{Alejandro Mahillo has been partly supported by Project ID2019-105979GBI00 of the MICINN of Spain and Project E48-20R, 
Gobierno de Aragón, Spain.\\ 
Silvia Rueda has been partly supported by Agencia Nacional de Investigación y Desarrollo (ANID), FONDECYT INICIACIÓN 2023, Grant 11230856.}

\subjclass[2020]{47A10, 47A35, 47D03}

\keywords{Ritt condition, Kreiss condition, spectral localization, norm estimates, Stolz domain}

\begin{abstract}
    A new condition is introduced by generalizing the Ritt and Kreiss operators named $(\al,\be)$-RK condition. Geometrical
    properties of the spectrum for the case $\be < 1$ are studied, moreover it is shown that in that case if $\al + \be = 1$ 
    the operator is Ritt. Estimates for the power and power differences norms for this type of operators are also studied. 
    Lastly we apply this theory to obtain and interpolation result over Ritt and Kreiss operator on $L^p$ spaces. 
\end{abstract}

\maketitle

\section{Introduction}
\label{sec:introduction}

The Kreiss matrix theorem holds an important place in the theory of finite difference methods for partial differential equations
as it deals with sufficient and necessary conditions about the stability of general systems. In particular, an upper bound is obtained
for the powers of an operator $T \in \BB(X)$, the space of bounded linear operators on the finite dimensional Banach space $X$. 
A bounded linear operator $T$ satisfies the \textit{Kreiss condition} if there exists $C \geq 1$ such that
\begin{equation}
    \label{eq:Kreiss_condition}
    \norm{(\lambda I - T)^{-1}} \leq \frac{C}{\abs{\lambda}-1}, \; \abs{\lambda}>1,
\end{equation}
where $\sigma(T) := \CC \setminus \rho(T)$ is the spectrum of $T$, $\rho(T)$ is the resolvent set and $\overline{\DD}$ stands
for the closure of the open unit disk $\DD$. When the dimension of the Banach space $X$ is finite the Kreiss matrix theorem states that 
\[
    \norm{T^n} \leq C e \dim X, \qquad n \in \NN,
\]
although, in general, in an infinite dimensional case a Kreiss operator satisfies $\norm{T^n} = O(n)$. See \cite{LyubichNevanlinna1991},
where it is proved and see \cite{Shields1978} where optimality is shown. The Kreiss matrix theorem  relates bounds of the resolvent
operator $R(\lambda, T) := (\lambda I - T)^{-1}$, with $\lambda \in \rho(T)$, to bounds of the norm of the powers of the operator.
Moreover, using the series expansion of the resolvent operator it can be proven that if an operator  $T$ is power bounded, that is,
$\sup_{n \in \NN} \norm{T^n} < \infty$, then $T$ satisfies the Kreiss condition.  See \cite{StrikwerdaWade1997} and references
therein for more information about Kreiss operators.

Another condition that appears in numerical analysis of the same kind is the Ritt condition. An operator $T$ is said to satisfy the
\textit{Ritt condition} if there exists $C \geq 1$ such that
\begin{equation}
    \label{eq:Ritt_condition}
  \norm{(\lambda I - T)^{-1}} \leq \frac{C}{\abs{\lambda-1}}, \; \abs{\lambda}>1.
\end{equation}
In \cite{Lyubich1999,NagyZemanek1999} is proven that if $T$ is Ritt, then it is power bounded and therefore it is Kreiss.
The Ritt condition can be reformulated in a more potent form: there exists an angle $\omega \in [0, \pi / {2})$ such that
\begin{equation}
    \label{eq:Ritt_condition_sectorial}
     \norm{R(\lambda,T)} \leq \frac{C_{\omega'}}{\abs{\lambda-1}}, \,
    \lambda \in \CC \setminus (1-\overline{\Sigma}_{\omega'}),
\end{equation}
for every $\omega' \in (\omega, \pi)$ and an appropriate $C_{\omega'} \geq 1$. Here, 
\[
    \Sigma_{\omega} = \{ \lambda \in \CC : \abs{\arg \lambda} < \omega\} \quad \text{and} \quad \Sigma_{0} = (0, \infty).
\]
In this situation, it is said that $T$ is of angle $\omega$. In particular, $\omega = \arccos(1 / {C})$, see \cite{Lyubich1999}.
This rephrasing of the conditions offers a more detailed understanding of the operator $T$ and its spectrum. Nagy and Zemánek 
\cite{NagyZemanek1999} and Lyubich \cite{Lyubich2001} established a characterization for Ritt operators: $T$ is Ritt if and only
if $\sup_{n \geq 1} \norm{T^n} < \infty$ and $\sup_{n \geq 1} n \norm{T^{n+1}-T^{n}} < \infty$. This implies that $T$
is power bounded with a decay rate in the Katznelson-Tzafriri theorem. However, power boundedness falls short of
satisfying the Ritt condition, as proved in \cite{NagyZemanek1999}. An operator may exhibit $\norm{T^{n+1}-T^{n}} \to 0$
without attaining power boundedness, as exemplified in \cite{Swiech1997}.

One could wonder if it is possible to have $\sup_{n \geq 1} n \norm{T^{n+1}-T^{n}} < \infty$ and $T$ not being power bounded.
As 
\[
    T^n = \frac{1}{n} \sum_{k=1}^n T^k + \frac{1}{n} \sum_{k=1}^{n-1} k (T^{k+1}-T^k)   
\]
we have that a mean bounded operator meeting $\sup_{n \geq 1} n \norm{T^{n+1}-T^{n}} < \infty$ is actually power bounded.
This question was debated in the early 2000s and definitively answered in the negative by
\cite[Theorem 3.3]{KaltonMontgomery-SmithOleszkiewiczTomilov2004}. There are many counterexamples, but the assumption on $T^{n+1}-T^{n}$
is extremal in the sense that $\liminf_n n \norm{T^{n+1}-T^{n}} > 0$ always holds (see \cite{Esterle1983}). There is a large body
of literature on Ritt operators, a sample of these might include the following titles \cite{BorovykhDrissiSpijker2000,ElFallahRansford2002,
GomilkoTomilov2018,LancienLeMerdy2015,Lyubich1999,NagyZemanek1999, Nevanlinna1993,Nevanlinna1997,Nevanlinna2001,
KaltonMontgomery-SmithOleszkiewiczTomilov2004,Dungey2011,Vitse2004}.

A Ritt operator is power bounded and therefore Kreiss. However, a Kreiss operator does not have to be Ritt since it does not 
even have to  be power bounded as mentioned above. This leads to the study of intermediate conditions between the Ritt 
and Kreiss conditions which could lead to new viewpoints in operator theory. We introduce the following new condition.

\begin{definition}
    \label{def:alpha_beta_RK_operator}
    An operator $T$ is said to be \textit{$(\al,\be)$-RK} with $\al,\be \geq 0$ if there exist a constant $C \geq 1$ such that,
    \begin{equation}
        \label{eq:alpha_beta_RK_operator}
        \norm{(\lambda I - T)^{-1}} \leq \frac{C \abs{\lambda}^{\al+\be - 1}}{\abs{\lambda-1}^\al (\abs{\lambda}-1)^\be},
        \; \abs{\lambda}>1.
    \end{equation}
\end{definition}

Borovykh and Spijker initially employed a more specific version of this condition in \cite{BorovykhSpijker2000Definition}.
Notably, they omitted the modulus of the value $\lambda$ in the numerator, treating it as bounded. Their study focused
solely on non-negative integer values of $\alpha$ and $\beta$ and operators on  finite dimensional Banach spaces. They also present
as an application of its research and example of an $(1,1)$-RK operator. In our work, we extend this definition to a more general context.
Vitse termed the less general version, the $(0,\alpha+1)$-RK condition, designating it as a generalized Kreiss
condition of order $\alpha$ in \cite[Theorem 2.6]{Vitse2005}.

Furthermore, Nevanlinna explored the $(\alpha+1,0)$-RK condition as an enhancement to a Kreiss operator in \cite[Theorem 9]{Nevanlinna2001}.
The subsequent extension to general bounded operators by Cohen and Lin in \cite{CohenLin2016} incorporated results from Seifert's
work in \cite[Corollary 3.1]{Seifert2016}. Seifert's research specifically addressed the rates of decay in the classical Katznelson-Tzafriri
theorem with a general decay of the resolvent on $\mathbb{T} \setminus \{1\}$.

A few remarks about the definition. The first one is about the constant, notice that we have chosen $C \geq 1$ and one may wonder why
not giving the definition for $C>0$. Notice that the operator $\lambda R(\lambda,T) \to I$ as $\abs{\lambda} \to \infty$ this force the
constant to be greater than or equal to 1. For the second remark we denote by $\mathcal{RK}_{\al,\be}$  the set of operators satisfying
the $(\al,\be)$-RK condition. Thus, if $0 \leq \al_0 \leq \al_1$ and $0 \leq \be_0 \leq \be_1$, we have that $\mathcal{RK}_{\al_0,\be_0}
\subset \mathcal{RK}_{\al_1,\be_1}$, and we also have due to the reverse triangle inequality that
\begin{equation}
    \label{eq:sets_RK}
    \mathcal{RK}_{\al + \be,0} \subset \mathcal{RK}_{\al ,\be} \subset \mathcal{RK}_{0,\al + \be}.
\end{equation}
In particular if $\al+\be=1$ the set $\mathcal{RK}_{\al ,\be}$ lies in between the Ritt operators and Kreiss operators.

The article focuses on two aspects of operators satisfying the $(\al,\be)$-RK condition: geometric properties of the spectrum
and growth of power norms. Lyubich showed that the spectrum of a Ritt operator lies within a sectorial region, while Gomilko and Tomilov
proved that it is within a Stolz domain $S_\sigma = \{z \in \DD \, : \, \abs{1-z}/(1-\abs{z}) < \sigma \} \cup \{1\}$. On the other hand,
we proceed as Borovykh and Spijker in \cite{BorovykhSpijker2000Definition} studying what happens to the power norms and the differences
of powers norm. Although it may seem that these two concepts are unrelated, obtaining information about the shape of the operator's
spectrum is fundamental for deriving accurate estimates of the operator's power norms. There are two key ways in which this relationship
is fruitful. The first is through Gelfand's formula, which establishes a connection between the spectral radius and the growth of an
operator's powers. The \textit{spectral radius} $r(T)$ of a bounded linear operator $T$ is defined as
$r(T) := \max \{ \abs{\lambda} : \lambda \in \sigma(T)\}$, where we recall that $\sigma(T)$ is a compact subset of $\CC$.
The formula is given by:
\begin{equation}
    \label{eq:Gelfands_formula}
    r(T) = \lim_{n \to \infty} \norm{T^n}^{1 / {n}}.
\end{equation}
The second method involves using the Cauchy integral formula to derive precise estimates of the norm of the operator's powers.
The key here is identifying the optimal curve for integration to yield the most accurate estimates.

The paper follows the following distribution. In \Cref{sec:geometrical_properties}, we extend some spectral results for Ritt 
operators to operators satisfying our new condition when $\be < 1$, aiming to define a functional calculus for these new operators.
Specifically, we show that the localization of the spectrum of an  $(\al,1-\al)$-RK operator for $\al \in (0,1]$ is analogous to the
localization of the spectrum of a Ritt operator. In \Cref{sec:norm_estimates}, we study the sequences $\{\norm{T^n}\}_{n\in \NN}$ and
$\{\norm{T^{n+1}-T^{n}}\}_{n \in \NN}$ for an operator $T$ that is $(\al,\be)$-RK. We find that if an operator is $(\al,1-\al)$-RK for
$\al \in (0,1]$, then it is Ritt, given that an operator is Ritt if and only if it is power bounded and
$\sup_{n \geq 0} n \norm{T^{n+1}-T^n} < \infty$ (see \cite[Theorem 4]{NagyZemanek1999}). \Cref{sec:applications} presents an 
application of our results to $L^p$ spaces.

\section{Geometrical properties of the spectrum}
\label{sec:geometrical_properties}

The aim of this section is to discuss geometrical properties of the spectrum of $(\al,\be)-RK$ operators. Lyubich (\cite{Lyubich1999})
proved that the spectrum $\sigma(T)$ of a Ritt operator is contained in a special convex domain. Moreover, he also proved that a
generalized Ritt condition leads to a similar localization result which cannot be applied directly to our new condition. Recently,
Gomilko and Tomilov (\cite{GomilkoTomilov2018}) proved a similar characterization of Ritt operators under the condition that the spectrum
is contained inside a Stolz domain $S_\sigma$. For $\sigma \geq 1$, the \textit{Stolz domain} $S_\sigma$ is defined as
\begin{equation}
    \label{eq:Stolz_domain}
    S_\sigma := \left\{ z \in \DD: \frac{\abs{1-z}}{1-\abs{z}} < \sigma \right\} \cup \{1\}.
\end{equation}
Clearly, $S_\sigma = \{1\}$ if $\sigma = 1$. As explained by them the choice of this kind of sets is beneficial in order to define a
proper functional calculus for Ritt operators. Thus, our approach would be similar to Gomilko and Tomilov's one.

In our case the spectrum of $(\al,\be)$-RK operators also have interesting geometrical properties when $\be < 1$. To do so, we will
start with a few elementary observations that follow from the properties of the resolvent of an operator. Recall that the resolvent set
$\rho(T)$ is an open set; moreover, if $\lambda \in \rho(T)$ then the right-hand side of
\begin{equation}
    \label{eq:resolvent_series}
    R(\mu,T) = \sum_{k=0}^\infty R(\lambda,T)^{k+1} (\lambda - \mu)^k
\end{equation}
is well-defined for $\mu$ inside the set 
\begin{equation}
    \label{eq:resolvent_set_inequality}
   \mathcal{S}:= \{ \mu : \abs{\mu - \lambda} < \norm{R(\lambda,T)}^{-1}\}.
\end{equation}

The right-hand side of \eqref{eq:resolvent_series} defines $[R(\lambda,T) (I - (\mu-\lambda) R(\lambda,T))]^{-1}$,
a product of invertible operators. It shows that $R(\mu,T)$ is well-defined and hence $\mu \in \rho(T).$ This shows
$\mathcal{S}\subset \rho(T)$. As a consequence, the distance from $\sigma(T)$ to a value in the resolvent set satisfies
\begin{equation}
    \label{eq:distance_spectrum}
    d(\lambda, \sigma(T) )\geq \norm{R(\lambda,T)}^{-1}, \quad \lambda\in  \rho(T).
\end{equation} 
It is important to mention that \eqref{eq:resolvent_series}, \eqref{eq:resolvent_set_inequality} and \eqref{eq:distance_spectrum}
were employed in \cite{Lyubich1999} to investigate the localization of the spectrum of Ritt operators. In our case, we have the 
following lemma.

\begin{lemma}
    Let $T$ be an $(\al,\be)$-RK operator with $\be<1$ then $\sigma(T) \subset \DD \cup \{1\}$.
\end{lemma}
\begin{proof} 
    Take $\lambda \in \TT \setminus \{1\}$ where $\TT:=\{\lambda \in \CC :\abs{\lambda}=1\}$ and suppose $\lambda \in \sigma(T)$.
    Moreover, consider the values $\lambda_{\e} = \left( 1+ \e \right) \lambda \in \rho(T)$ with $\e>0$. Thus, from 
    \eqref{eq:distance_spectrum} we have that 
    \[
        \e \geq \norm{R(\lambda_{\e},T)}^{-1} \geq C^{-1} \abs{\lambda_{\e}-1}^\al \e^\be \left(1+\e\right)^{1-\al-\be}.
    \]
As $\e \to 0$ we get a contradiction. Thus, $\lambda \in \rho(T)$.
\end{proof}

Moreover, we have the following result.

\begin{proposition}
    \label{pr:RK_resolvent_estimate}
    Let $T$ be an $(\al,\be)$-RK operator with $\be<1$, then we have the following resolvent estimates,
    \begin{align}
        \label{eq:bounds_on_the_torus}
            \norm{R(\lambda, T)} & \leq \frac{C_{\al,\be} }{ \abs{\lambda-1}^\frac{\al}{1-\be}}, \qquad \lambda \in \TT \setminus \{1\},         
        \end{align}
    where, we have the following cases:
    \begin{enumerate}
        \item If $\be=0$ then 
        \[
            C_{\al,\be} = C   
        \]
        i.e. $C_{\al,0} = C$ for every $\alpha$.
        \item If $\be>0$ and $\al + \be \leq 1$ then 
                \[
                    C_{\al,\be} = \frac{C^\frac{1}{1-\be}}{\be^{\frac{\be}{1-\be}}(1-\be)}.
                \]
        \item If $\be>0$ and $\al + \be > 1$ then
                \[
                    C_{\al,\be} = \frac{C^\frac{1}{1-\be}\, 2^{\frac{\al}{1-\be}}}{\be^{\frac{\be}{1-\be}}(1-\be)}.
                \]
    \end{enumerate}
\end{proposition}
\begin{proof}
    We start the proof by taking $\lambda \in \TT \setminus \{1\}$ and the values
    $\lambda_{\e} = \left( 1+ \e \right) \lambda \in \rho(T)$ with $\e>0$. The resolvent identity seen in \eqref{eq:resolvent_series}
    together with the condition of being an $(\al,\be)$-RK operator yields,
    \begin{equation}
        \label{eq:general_resolvent_series}
        \begin{aligned}
            \norm{R(\lambda,T)} &   \leq    \sum_{k=0}^\infty \abs{\lambda - \lambda_\e}^k \norm{R(\lambda_\e,T)}^{k+1} \\
                                &   \leq    \frac{C(1+\e)^{\al+\be-1}}{\abs{\lambda_\e-1}^\al \e^\be} 
                                            \sum_{k=0}^\infty \left(\frac{C(1+\e)^{\al+\be-1}\e^{1-\be}}{\abs{\lambda_\e-1} ^\al}\right)^k \\
                                &   =       \frac{(1+\e)^{\al+\be+1}}{q \abs{\lambda_\e-1}^\al \e^\be - (1+\e)^{\al+\be+1} \e},
        \end{aligned}
    \end{equation}
    where $q = 1 / {C}$. Notice that we can choose $\e$ small enough so that series converges. Now we split the study of this
    value in different cases.

    \textit{Case $\be = 0$:} In this case, \eqref{eq:general_resolvent_series} becomes
    \[
        \norm{R(\lambda,T)} \leq    \frac{(1+\e)^{\al+1}}{q \abs{\lambda_\e-1}^\al - (1+\e)^{\al+1} \e}   
    \]
    and taking $\e \to 0$ gives,
    \[
        \norm{R(\lambda,T)} \leq \frac{C}{\abs{\lambda-1}^\al}.
    \]

    \textit{Case $\be > 0$ and $\al + \be \leq 1$:} Now, by a geometrical reasoning (you can draw the vectors), we note that 
    $\abs{\lambda - 1} \leq \abs{\lambda_\e - 1}$ and as $\al + \be \leq 1$ we have that $(1+\e)^{\alpha+\beta-1} \leq 1$.
    Thus, \eqref{eq:general_resolvent_series} becomes
    \[
        \norm{R(\lambda,T)} \leq \frac{1}{q \abs{\lambda-1}^\al \e^\be - \e},  
    \]
    but in this case we can't use a limit argument as before. In this case we find the best value $\e$ which minimizes the
    right-hand side of the last inequality. Define the function 
    \[
        f_\lambda(\e) := \left( q \abs{\lambda-1}^\al \e^\be - \e \right)^{-1}
    \]
    and remember that $\e \in (0,r)$, for some $r>0$. We compute the derivative of $f_\lambda$
    \[
        f_\lambda'(\e)= -   \left(q \abs{\lambda-1}^\al \e^\be - \e\right)^{-2}
                            \left(q \be \abs{\lambda-1}^\al \e^{\be-1} - 1\right).
    \]
    Thus, the minimum is reached at $\e^* = \left(q \be \abs{\lambda-1}^\al \right)^{1/{(1-\be)}}$ and notice that for this value
    the series appearing in \eqref{eq:general_resolvent_series} is convergent. Finally, we compute $f_\lambda(\e^*)$ to find the
    desired bound,
    \[
        f_\lambda(\e^*) =   \left( \be^\frac{\be}{1-\be} (q \abs{\lambda-1}^\al)^\frac{1}{1-\be} 
                            - \left(q \be \abs{\lambda-1}^\al \right)^{1/{(1-\be)}} \right)^{-1}
                        =   \frac{C^\frac{1}{1-\be}}{\be^\frac{\be}{1-\be} (1-\be) \abs{\lambda-1}^\frac{\al}{1-\be}}.
    \]

    \textit{Case $\be > 0$ and $\al + \be > 1$:}
    The ideas in this case are similar as the one that we just did. The problem here is that if we use the same techniques as before
    the optimization problem that we find is way more complicated to find an exact solution, so in order to solve that problem
    we prove the following inequality
    \begin{equation}
        \label{eq:ineq_complex_numbers}
        \frac{1+\e}{\abs{\lambda_\e -1}} \leq \frac{2}{\abs{\lambda-1}}. 
    \end{equation}
    Indeed, by the reverse triangle inequality we have that $\abs{\lambda_\e -1} \geq \e$, therefore,
    \[
        \frac{(1+\e) \abs{\lambda-1}}{\abs{\lambda_\e-1}}   =       \frac{\abs{\lambda_\e-1-\e}}{\abs{\lambda_\e-1}}
                                                            \leq    1 + \frac{\e}{\abs{\lambda_\e-1}} \leq 2,
    \]
    and the result follows. Thus, taking into account that $(1+\e)^{\be-1}<1$ and inequality \eqref{eq:ineq_complex_numbers},
    the bound \eqref{eq:general_resolvent_series} becomes,
    \[
        \norm{R(\lambda,T)} \leq \frac{1}{2^{-\al} q \abs{\lambda-1}^\al \e^\be - \e}.  
    \]
    Notice that the optimization problem is, indeed, the same as in the last case. We have just changed the value $q$ 
    by $q 2^{-\alpha}$. Therefore, the minimum is reached at
    $\e^* = \left(q 2^{-\alpha} \be \abs{\lambda-1}^\al \right)^{1/{(1-\be)}}$, and
    \[
        \norm{R(\lambda,T)} \leq \frac{C^\frac{1}{1-\be} 2^\frac{\alpha}{1-\be}}
                                    {\be^\frac{\be}{1-\be} (1-\be) \abs{\lambda-1}^\frac{\al}{1-\be}}.
    \] 
\end{proof}

This result provides important insights about the spectrum's shape as we approach $1$. In \cite[Theorem 3]{Lyubich1999}, 
a general method is proposed for operators that meet certain conditions. Specifically, these operators must satisfy bounds like
$\norm{R(\lambda,T)} \leq C / {\gamma(\theta)}$, where $\lambda = e^{i \theta}$, $\theta \in (0,2\pi)$, and $\gamma$ is a function that
belongs to $C^2[0,2\pi]$, the set of functions with continuous second derivatives on the interval $[0,2\pi]$. The bounds in equation
\eqref{eq:bounds_on_the_torus} can be seen as a special case of this condition. Here, $\abs{\lambda-1}^{\al / {(1-\beta)}} 
= 2^{\al / {(1-\beta)}}\sin(\theta / {2})^{\al / {(1-\beta)}}$ with $\lambda = e^{i \theta}$ and $\theta \in (0,2\pi)$. 
However, this theorem can't be directly applied to our case. The reason is that it requires the function $\gamma$ to belong to the 
space $C^2[0,2\pi]$,  which might not always be true in our situation. 

\begin{figure}[htb]
    \centering
    \includegraphics[width=0.7\textwidth] {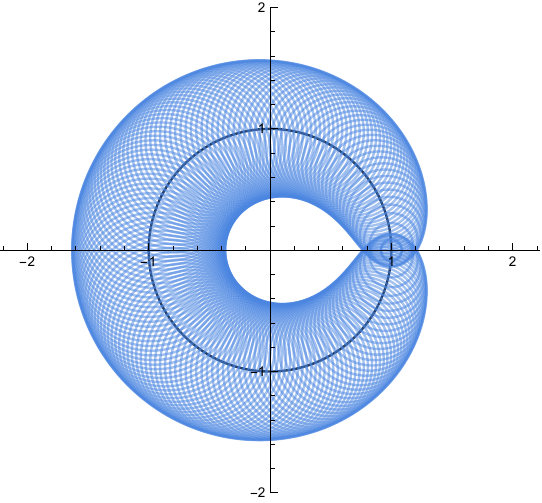}
    \caption{Case $\alpha=1 / 4$, $\beta=1 / 4$ and $C=2$. The spectrum must lie inside the white region inside the unit disk.}
    \label{fig:DecayCase}
\end{figure}

Now we introduce the special case when $\al + \be < 1$ which follows directly from \Cref{pr:RK_resolvent_estimate}. To give an idea of
the proof you can check \Cref{fig:DecayCase} where in each point of the torus we have drawn the circumference given by the estimate
computed in \Cref{pr:RK_resolvent_estimate}. Therefore, the spectrum has to be outside the blue circles and inside the unit disk.

\begin{proposition}
    Suppose $T$ is an $(\al,\be)$-RK operator where $\be<1$ and $\al + \be < 1$. Then, we have $\sigma(T) \subset \mathbb{D}$.
\end{proposition}
\begin{proof}
    Based on \Cref{pr:RK_resolvent_estimate}, we know that $\sigma(T) \subset \mathbb{D} \cup \{1\}$. Now, let's assume
    $1 \in \sigma(T)$ and take a sequence $\{\lambda_n\}_{n \in \NN} \subset \mathbb{T} \setminus \{1\}$ such that
    $\lim_{n \to \infty} \lambda_n = 1$. From property \eqref{eq:resolvent_set_inequality}, using \eqref{eq:bounds_on_the_torus}
    we get 
    \[
        |\lambda_n-1|^{-1} \leq \|R(\lambda_n,T)\| \leq C |\lambda_n-1|^{-\frac{\al}{1-\be}}.
    \]
    This implies,
    \[
        C \geq |\lambda_n-1|^{\frac{\al+\be-1}{1-\be}},
    \]
    leading to a contradiction since $|\lambda_n-1|^{\frac{\al+\be-1}{1-\be}} \to \infty$ as $n \to \infty$ (given that $\al + \be < 1$).
\end{proof}

\begin{remark}
    \label{re:exp_decay}
    This result leads us to the conclusion that an $(\al,\be)$-RK operator, given $\be<1$ and $\al + \be < 1$, verifies that $r(T)<1$.
    According to Gelfand's formula \eqref{eq:Gelfands_formula}, this means the sequence $\{\norm{T^n}\}_{n \in \NN}$ decays exponentially.
\end{remark}

\begin{figure}[htb]
    \centering
    \includegraphics[width=0.7\textwidth] {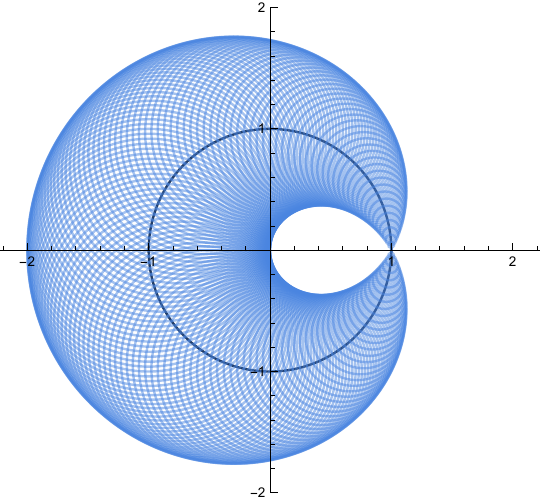}
    \caption{Case $\alpha=1 / 2$, $\beta=1 / 2$ and $C=2$. The spectrum must lie inside the white region inside the unit disk.}
    \label{fig:RittCase}
\end{figure}

\begin{remark}
    \label{re:ritt_geom}
    From \Cref{pr:RK_resolvent_estimate} we also have that if $\al + \be = 1$, with $\be<1$ then our operator satisfies the Ritt
    condition on $\TT$ and therefore as in \cite[Proposition 1]{Lyubich1999} we have that $\sigma(T) \subset (1-\overline{\Sigma}_{\omega})$
    with $\omega = \arccos(1 / {C_{\al,\be}})$. Or if we use the notion of Stolz domain \eqref{eq:Stolz_domain} we have that 
    $\sigma(T) \subset \overline{S}_{C_{\al,\be}}$. See \Cref{fig:RittCase} for a visual example. The spectrum is contained in the region
    inside the unit disk painted in white.
\end{remark}

\begin{figure}[htb]
    \centering
    \includegraphics[width=0.7\textwidth] {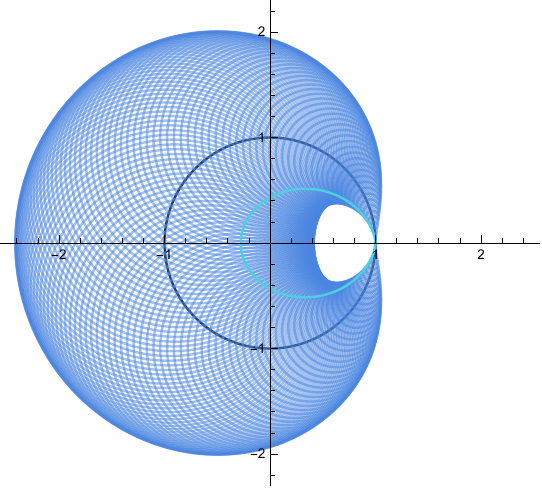}
    \caption{Case $\alpha=3 / 4$, $\beta=1 / 2$ and $C=2$, in cyan is colored the border of the $\frac{\al}{1-\beta}$-Stolz
    region. The spectrum of $T$ is contained in the white region inside $\DD$.}
    \label{fig:StolzCase}
\end{figure}

For the moment, we have described the geometrical shape of the spectrum for the cases when $\be<1$ and $\al + \be \leq 1$. 
The next step will be to present what happens for the case $\be<1$ and $\al + \be > 1$. The idea is to extend the definition of Stolz
domain \eqref{eq:Stolz_domain}, for which we introduce the concept of $\al$-Stolz domains. For $\sigma > 0$ we define the $\al$-Stolz 
domain $S_\sigma^\al$ by
\[
    S_\sigma^\al := \left\{ z \in \DD: \frac{\abs{1-z}^\al}{1-\abs{z}} < \sigma \right\} \cup \{1\}.
\]
Notice that $S_\sigma^\al$ is a convex set for $\al \geq 1$, and if $\sigma > 1$ then $0 \in S_\sigma^\al$.
See \Cref{fig:StolzCase} for a visual example. In this context we have the following theorem.

\begin{theorem}
    Let $T$ be an $(\al,\be)$-RK operator with $\be<1$ and $\al + \be >1$ then,
    \begin{equation}
        \label{eq:geometric_condition_Stolz_Ritt}
        \sigma(T)  \subset \overline{S}_\sigma^\frac{\al}{1-\be}, \qquad \text{with } \sigma = 2 C_{\al,\be}, \\
    \end{equation}
    with $C_{\al,\be}$ defined as in \Cref{pr:RK_resolvent_estimate}.
\end{theorem}
\begin{proof}
    Our first step is to find a lower bound for the constant $C_{\al,\be}$. In order to have a non-empty spectrum the constant
    $C_{\al,\be}$ must verify that 
    \begin{equation}
        \label{eq:constant_lower_bound}
        C_{\al,\be} \geq 2^{\frac{\al}{1-\be}-1}.   
    \end{equation}
    Indeed, if $C_{\al,\be} < 2^{\frac{\al}{1-\be}-1}$ by \Cref{pr:RK_resolvent_estimate} we have that
    $\norm{R(-1,T)} \leq C_{\al,\be} 2^{-\frac{\al}{1-\be}}$ and \eqref{eq:distance_spectrum} yields,
    \[
        d(-1,\sigma(T)) \geq \norm{R(-1,T)}^{-1} \geq C_{\al,\be}^{-1} 2^{\frac{\al}{1-\be}} > 2,
    \]
    where we have used the assumption that $C_{\al,\be} < 2^{\frac{\al}{1-\be}-1}$. This inequality implies that
    $\bar{\DD} \subset \rho(T)$ which implies that the spectrum is empty, which is a contradiction.

    We prove the case $\be = 0$ as the ideas are the same, but the text is going to be clearer for the reader. 
    For the rest of the proof we denote by $q = 1 /{C}$. We define the following set,
    \[
        \Omega(q,\al)= \{ r e^{i\theta} \in \DD : r \leq 1 - q \abs{e^{i \theta}-1}^\al \} \cup \{1\}.
    \]
    Notice that $\sigma(T) \subset \Omega(q,\al)$. From \Cref{pr:RK_resolvent_estimate} together with
    inequality \eqref{eq:distance_spectrum} applied to the points $e^{i\theta} \in \rho(T)$ and
    $1 \neq r e^{i\theta} \in \sigma(T)$ yields
    \[
        1-r = d(e^{i\theta}, re^{i\theta}) \geq d(e^{i\theta}, \sigma(T)) \geq \norm{R(e^{i\theta},T)}^{-1} \geq q \abs{e^{i \theta}-1}^\al. 
    \]
    Thus, $r e^{i\theta} \in \Omega(q,\al)$.

    The final step is to prove that $\Omega(q,\al) \subset \overline{S}_\sigma^\al$. In order to achieve this notice that we just
    have to prove it for the border of the set $\Omega(q,\al)$ as we know that $\overline{S}_\sigma^\al$ is convex and
    $ 0 \in \overline{S}_\sigma^\al$. Moreover, it's enough to prove it just for the complex value with positive imaginary part as
    both sets are symmetrical about the real axis. Thus, take  $z = r(\theta) e^{i \theta} \in \Omega(q,\al)$ with
    $r(\theta) = 1-q \abs{e^{i \theta}-1}^\al$ for $\theta \in (0,\pi)$. Therefore,
    \[
        \begin{aligned}
            \frac{\abs{1-z}^2}{(1-\abs{z})^{2/ {\al}}}  & = \frac{(1-r(\theta))^2 + 4r(\theta)\sin(\theta/{2})^2}{(1-r(\theta))^{2/ {\al}}}\\
                                                        & = \frac{1}{q^{2/ {\al}}} \left( 1 + q^2 4^{\al-1} \sin(\theta / {2})^{2\al-2}
                                                            - q 2^\al \sin(\theta / {2})^\al \right),
        \end{aligned}
    \]
    where we have used the fact that $\abs{e^{i \theta}-1}^\al = 2^\alpha \sin(\theta / {2})^\al$. In order to make it easier
    to handle we substitute $x = \sin(\theta / {2}) \in (0,1)$, and we study the maximum value of the function,
    \[
        f(x) = 1 + q^2 4^{\al-1} x^{2\al-2} -q 2^\al x^\al,
    \]
    with derivative,
    \[
        f'(x)  = q 2^\al (q(\al-1)2^{\al-1}x^{\al-2}-\al)x^{\al-1},
    \]
    in the interval $(0,1)$. If $\al \geq 2$ then $q(\al-1)2^{\al-1}x^{\al-2}-\al < 0$ so f is decreasing in $(0,1)$ and
    \[
        \frac{\abs{1-z}^2}{(1-\abs{z})^{2/ {\al}}} < \frac{1}{q^{2/ {\al}}}.
    \]
    If $\al \in (1,2)$ we have a maximum at $x_{max} = \left( \frac{q(\al-1)2^{\al-1}}{\al}\right)^{1/ {(2-\al)}}$.
    Notice that $x_{max} \in (0,1)$. We compute $f(x_{max})$ to obtain that,
    \[
        f(x_{max}) = 1 + q^\frac{2}{2-\al} \left(\frac{2(\al-1)}{\al}^{\frac{2\al-2}{2-\al}}\right)
                    \left(1-\frac{2(\al-1)}{\al}\right),
    \] 
    as $q \leq 2^{1-\al}$ by \eqref{eq:constant_lower_bound}, $f(x_{max}) \leq 2^{2/ {\al}}$ with $\al \in (1,2)$.
    
    To sum up, we have that, 
    \[
        \frac{\abs{1-z}^2}{(1-\abs{z})^{2/ {\al}}} \leq \frac{2^{2/ {\al}}}{q^{2/ {\al}}},
    \]
    and the claim follows.
\end{proof}

\begin{remark}
    In \cite{CohenLin2016} they study a similar condition to ours but without the Kreiss component in the bound. In order to
    localize the spectrum of the operators they introduced a set called \it{quasi-Stolz} set. Their starting point is not our 
    Stolz domain because, for them, the geometric construction of the Stolz region is by taking a circle of radius $r<1$ 
    centered at $0$ and drawing two tangent line segments from the point $1$ to this circle. Their generalization comes from
    the construction of Paulauskas \cite{Paulauskas2012} where the tangent lines are replaced by arcs of a \it{tangent}
    ``parabola-like'' curve $x=1-b\abs{y}^\alpha$, $1<\alpha<2$, $b>0$, or $\alpha = 2$ and $b>\frac{1}{2}$ (with 
    $\abs{y} \leq \abs{y_0} < 1$). The convex set obtained is what they call a quasi-Stolz set.
\end{remark}

\section{Norm estimates of powers of \texorpdfstring{$(\al,\be)$}{(α, β)}-RK operators}
\label{sec:norm_estimates}

In this section, we aim to understand how the new condition on the resolvent influences the boundaries we can set on the norm 
of the powers of our operator and other estimates. Before we delve into the norm estimates of such operators, we introduce two
integral estimates that will play a pivotal role in our subsequent findings.

\begin{lemma}
    \label{lemma:integral_upper_bound}
    Define the following integrals,
    \[
        I_\gamma(r) = \int_{-\pi}^\pi \frac{1}{\abs{re^{it}-1}^\gamma} \, dt, \; \; 
        J_\gamma(r) = \int_{-\pi}^\pi \frac{\abs{e^{it}-1}}{\abs{re^{it}-1}^\gamma} \, dt, \quad
        \text{for } r \in (1,1 + 1 / {2}).
    \]
    Then we have the following results:
    \begin{enumerate}
        \item There exists a constant $M_\gamma>0$ depending on $\gamma$ such that,
        \begin{enumerate}
            \item $I_\gamma(r) \leq M_\gamma$, for $0 \leq \gamma< 1 $,
            \item $I_1 (r) \leq M_1 \log((r-1)^{-1})$ when $\gamma = 1$,
            \item $I_\gamma(r) \leq M_\gamma (r-1)^{-\gamma+1}$ when $\gamma > 1$.
        \end{enumerate}
        \item There exists a constant $N_\gamma>0$ depending on $\gamma$ such that,
        \begin{enumerate}
            \item $J_\gamma (r) \leq N_\gamma$, for $ 0 \leq \gamma < 2$,
            \item $J_1 (r) \leq N_1 \log((r-1)^{-1})$ when $\gamma = 2$,
            \item $J_\gamma (r) \leq N_\gamma (r-1)^{-\gamma+1}$ when $\gamma > 2$. 
        \end{enumerate}
    \end{enumerate}
\end{lemma}
\begin{proof}
    We start computing the bounds for $I_\gamma(r)$. Using the parity of the function we have that,
    \[
        \int_{-\pi}^\pi \frac{1}{\abs{re^{it}-1}^\gamma} \, dt  
            = 2 \int_0^\pi \frac{1}{(r^2 + 1 - 2r\cos t)^\frac{\gamma}{2}} \, dt
            = 2 (r - 1)^{-\gamma} \int_0^\pi \frac{1}{(1 + \frac{2r(1-\cos t)}{(r - 1)^2})^\frac{\gamma}{2}} \, dt.
    \]
    As know that for all the cases $\gamma \geq 0$ we have on one hand, that $2(1- \cos t) \geq \frac{4t^2}{\pi^2}$ for
    $ 0 \leq t \leq \pi$, and on the other hand, $r>1$, thus
    \[
        I_\gamma(r) \leq 2 (r - 1)^{-\gamma} \int_0^\pi \frac{1}{(1 + \frac{4 t^2}{\pi^2 (r-1)^2})^\frac{\gamma}{2}} \, dt
                 =  \pi (r - 1)^{-\gamma + 1} \int_0^\frac{2}{r-1} \frac{1}{(1 + u^2)^\frac{\gamma}{2}} \, du.
    \]
    Now, we split the proof in 3 cases. First, take $0 \leq \gamma< 1 $ and use the fact that $1+u^2 \geq u^2$, therefore
    \[
        I_\gamma(r) \leq \pi (r - 1)^{-\gamma + 1} \int_0^\frac{2}{r-1} \frac{1}{u^\gamma} \, du
                    = \frac{ 2^{1-\gamma} \pi}{1-\gamma}, \quad \gamma \in [0,1).
    \]
    If $\gamma = 1$,
    \[
        I_1(r) \leq \pi \int_0^\frac{2}{r-1} \frac{1}{(1 + u^2)^\frac{1}{2}} \, du 
            = \pi \log\left(\frac{2}{r-1} + \sqrt{\frac{4}{(r-1)^2}+1}\right) 
    \]
    as $r<1+\frac{1}{2}$ we have that there exists a constant such that,
    \[
        I_1(r) \leq M_1 \log((r-1)^{-1}).
    \]
    For the last case $\gamma > 1$ we bound the integral in an larger interval,
    \[
        I_\gamma(r)   \leq \pi (r - 1)^{-\gamma + 1} \int_0^\infty \frac{1}{(1 + u^2)^\frac{\gamma}{2}} \, du
                    \leq \pi (r - 1)^{-\gamma + 1} \left( 1 + \int_1^\infty \frac{1}{u^\gamma} \, du \right), 
    \]
    thus there exists a constant $M_\gamma>0$ such that,
    \[
        I_\gamma(r)\leq M_\gamma (r - 1)^{-\gamma + 1}, \quad \gamma > 1.
    \]

    We now proceed for the bound for $J_\gamma(r)$. Following the same steps as in the other case and taking into account that
    $2(1- \cos t) \leq t^2$ for $ 0 \leq t \leq \pi$ we have,
    \[
        J_\gamma(r) = 2 (r - 1)^{-\gamma} 
                    \int_0^\pi \frac{(2(1-\cos t))^\frac{1}{2}}{(1 + \frac{2r(1- \cos t)}{(r - 1)^2})^\frac{\gamma}{2}}\, dt
                 \leq 2 (r - 1)^{-\gamma} \int_0^\pi \frac{t}{(1 + \frac{4t^2}{\pi^2 (r - 1)^2})^\frac{\gamma}{2}}\, dt.
    \]
    Now a change of variable gives,
    \[
        J_\gamma(r) \leq 
            \frac{\pi^2}{4} (r - 1)^{-\gamma+2} \int_0^{\left(\frac{2}{r-1}\right)^2} \frac{1}{(1 + u)^\frac{\gamma}{2}}\, dt.
    \]
    From here we can use the same techniques as we have used for the bounds of the integral $I_\gamma(r)$, and we will get the
    desired bounds.
\end{proof}

\begin{theorem}
    \label{th:powers_growth}
    Let $T$ be an $(\al,\be)$-RK operator and $k \in \NN \cup \{0\}$. Then, there is $M_{\al,\be,k} >0$ depending on 
    $\al,\be$ and $k$ such that for all $n \in \NN$ sufficiently big,
    \begin{itemize}
        \item if $\al(k+1) < 1$ then $\norm{T^n} \leq M_{\al,\be,k} \, n^{k (\be-1)+\be}$,
        \item if $\al(k+1) = 1$ then $\norm{T^n} \leq M_{\al,\be,k} \, n^{k (\be-1)+\be} \log(n)$,
        \item if $\al(k+1) > 1$ then $\norm{T^n} \leq M_{\al,\be,k} \, n^{(\al + \be-1)(k+1)}$.
    \end{itemize}
\end{theorem}
\begin{proof}
    
    The idea of this proof is to use the Cauchy's differentiation formula in order to represent the powers of an operator.
    Similarly to \cite{BorovykhDrissiSpijker2000}, after $k$-times integration by parts we have
    \[
        T^n = \binom{n+k}{k}^{-1} \frac{1}{2 \pi i} \int_{\Gamma} \lambda^{n+k} R(\lambda,T)^{k+1} \, d\lambda, \quad n \in \NN,
    \]
    where $\Gamma = \{ \lambda \in \CC : \abs{\lambda} = r\}$ for some $r>1$. In particular, we will take $r = 1 + \frac{1}{n}$.
    Applying now the inequalities of the norm and the resolvent condition we have that,
    \[
        \norm{T^n} \leq \binom{n+k}{k}^{-1} \frac{1}{2 \pi}
                        \int_{-\pi}^\pi \frac{C^{k+1} r^{n+k+(\al+\be-1)(k+1)}}{\abs{r e^{it} - 1}^{\al (k+1)}(r-1)^{\be (k+1)}} \, dt,
    \]
    and in particular using the definition of $I_{\al(k+1)}(r)$ seen in \Cref{lemma:integral_upper_bound} we have that,
    \begin{equation}
        \label{eq:bound_powers_intermediate_formula}
        \norm{T^n} \leq \binom{n+k}{k}^{-1} \frac{C^{k+1} r^{n+(\al+\be)(k+1)}}{2 \pi (r-1)^{\be (k+1)}} \, I_{\al(k+1)}(r).
    \end{equation}
    Lastly we bound the terms appearing in the right-hand side of the last inequality. First, as seen in 
    \cite[Equation 1.18, page 77, Vol.1]{Zygmund1968} we have that
    \[
        \binom{n+k}{k} = \frac{n^k}{k!}(1+ O(n^{-1}))  
    \]
    so, we have that $\binom{n+k}{k}^{-1} \leq A_k n^{-k}$ for some positive constant $A_k$. As $r = 1 + \frac{1}{n}$ 
    we have that $r^n = (1+\frac{1}{n})^n < e$. Also related with the value $r$ we have that 
    $r^{(\al+\be-1)(k+1)} = (1+\frac{1}{n})^{(\al+\be)(k+1)}$ is bounded by some constant as $k$ is fixed and $n$ is 
    going to be big. Thus, \eqref{eq:bound_powers_intermediate_formula} transforms into,
    \[
        \norm{T^n} \leq M_{\al,\be,k} n^{\be(k+1)-k} I_{\al(k+1)}(1+n^{-1}),
    \]
    where in $M_{\al,\be,k}$ are included all the constant that we have mentioned in the last paragraph. Lastly, we apply 
    the bounds seen in \Cref{lemma:integral_upper_bound}.
    \begin{itemize}
        \item If $\al(k+1) < 1$, then 
            \[
                \norm{T^n} \leq M_{\al,\be,k} \, n^{\be(k+1)-k} = M_{\al,\be,k} \, n^{(\be-1)k+\be}.
            \]
        \item If $\al(k+1) = 1$, then 
            \[
                \norm{T^n} \leq M_{\al,\be,k} \, n^{\be(k+1)-k} \log(n) = M_{\al,\be,k} \, n^{(\be-1) k+\be} \log(n).
            \]
        \item If $\al(k+1) > 1$, then 
            \[
                \norm{T^n} \leq M_{\al,\be,k} \, n^{\be(k+1)-k} n^{\al(k+1)-1} = M_{\al,\be,k} \, n^{(\al + \be-1)(k+1)}.
            \]
    \end{itemize}
    Where the constants of the \Cref{lemma:integral_upper_bound} enter into the constant $M_{\al,\be,k}$.
\end{proof}

As a corollary we have the following result to summarize some of the most common cases.

\begin{figure}[htb]
    \begin{minipage}{.62\textwidth}
        \vspace*{1cm}
        \centering
        \captionsetup{type=figure}
        \includegraphics[width=0.85\textwidth]{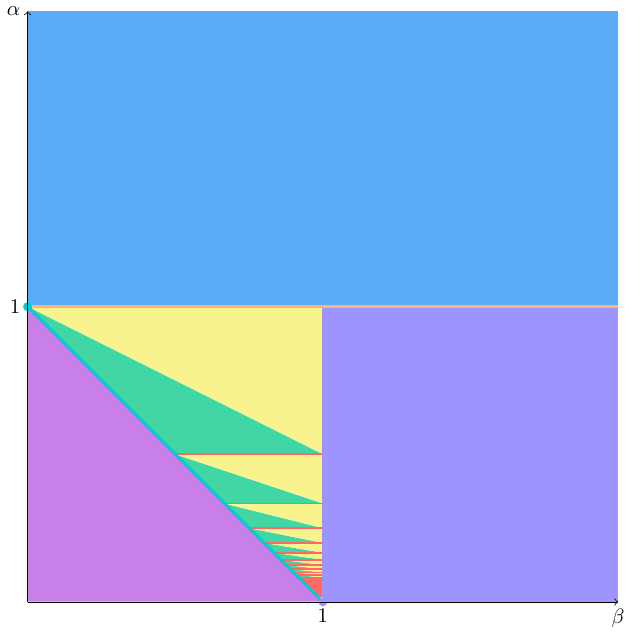}
        \captionof{figure}{Regions of \Cref{cor:regions_norms}.}
        \label{fig:regions_operator_powers}
    \end{minipage}
    \begin{minipage}{.33\textwidth}
        \begin{itemize}
            \item[\textcolor{newpurple}{\textbullet}] Exponential decay
            \item[\textcolor{newcyan}{\textbullet}] $O(1)$
            \item[\textcolor{newred}{\textbullet}] $O(n^{\frac{\alpha+\beta-1}{\alpha}} \log(n))$
            \item[\textcolor{newgreen}{\textbullet}] $O(n^{\floor{\frac{1}{\alpha}}(\beta-1)})$
            \item[\textcolor{newyellow}{\textbullet}] $O(n^{(\floor{\frac{1}{\alpha}}+1)(\alpha+\beta-1)})$
            \item[\textcolor{newviolet}{\textbullet}] $O(n^\beta)$
            \item[\textcolor{neworange}{\textbullet}] $O(n^\beta \log(n))$
            \item[\textcolor{newblue}{\textbullet}] $O(n^{\alpha+\beta-1})$
        \end{itemize}
    \end{minipage}
    \vspace*{1cm}
\end{figure}

\begin{corollary}
    \label{cor:regions_norms}
    Let $T$ be a $(\al,\be)$-RK operator on a complex Banach space $X$ with $\al,\be \geq 0$.
    Then, the norm of the powers of $T$ have the following growth.
    \begin{itemize}
        \item Case 1: $0 \leq \al < 1$ and $\al + \be < 1$. Then, there exists $\omega>0$ such that $\norm{T^n} = O(e^{-\omega n})$.
        \item Case 2: $\al = 0$ and $\be \geq 1$. Then, $\norm{T^n} = O(n^\be)$. In this case for $\beta=1$ we obtain the known
              growth for Kreiss operators.
        \item Case 3: $0 < \al < 1$ and $\al + \be = 1$. Then, $\norm{T^n} = O(1)$.
        \item Case 4: $0 < \al < 1$, $0 < \be < 1$  and $\al + \be > 1$.
        \begin{itemize}
            \item Case 4.1: $\al = \frac{1}{m+1}$ with $m \in \NN$. Then $\norm{T^n} = O(n^{\frac{\al+\be-1}{\al}}\log(n))$.
            \item Case 4.2: $\al$ is not the inverse of an integer.
            \begin{itemize}
                \item Case 4.2.1: $\frac{1}{\al}-\floor{\frac{1}{\al}} \geq \frac{\al+\be-1}{\al}$
                        then $\norm{T^n} = O(n^{(\al+\be-1)\left(\floor{\frac{1}{\al}}+1\right)})$.
                \item Case 4.2.2: $\frac{1}{\al}-\floor{\frac{1}{\al}} \leq \frac{\al+\be-1}{\al}$
                        then $\norm{T^n} = O(n^{\floor{\frac{1}{\al}}(\be-1)+1})$.
            \end{itemize}
        \end{itemize}
        \item Case 5: $0 < \al < 1$ and $\be \geq 1$. Then, $\norm{T^n} = O(n^\be)$.
        \item Case 6: $\al = 1$ and  $\be = 0$. Then, $\norm{T^n} = O(1)$. This case is known power boundedness of Ritt operators.
        \item Case 7: $\al = 1$ and  $\be > 0$. Then, $\norm{T^n} = O(n^\be \log(n))$.
        \item Case 8: $\al > 1$ and $\be \geq 0$. Then, $\norm{T^n} = O(n^{\al+\be-1})$.
    \end{itemize}
\end{corollary}
\begin{proof}
    The idea of the proof is to find the best integer $k$ to obtain the best growth in \Cref{th:powers_growth}.
    \begin{itemize}
        \item Case 1: $0 \leq \al < 1$ and $\al + \be < 1$. For this case we recall the following fact.
        If for some $N \geq 1 $, $\norm{T^N} < 1$, then by representing  $n$ as $n = sN + m$ with $0 \leq m < N-1$ we have that
        \[
            \norm{T^n} \leq \max_{0^\leq m < N} \norm{T^m} \norm{T^N}^s,
        \]
        and we conclude the exponential decay. We will prove now the case $\al = 0$ the case $\al>0$ follows the same idea.
        In this case $\norm{T^n} \leq M_{\al,\be,k} n^{k(\be-1)+\be}$ for every $k \in \NN \cup \{0\}$. As $\be<1$ we can take
        $k$ big enough so $k(\be-1)+\be<0$. Therefore, it exists $N \geq 1 $, such that $\norm{T^N} < 1$ and the 
        exponential decay follows.

        \item Case 2: $\al = 0$ and $\be \geq 1$. The best growth is achieved for $k=0$, and we have that $\norm{T^n} = O(n^{\be})$.

        \item Case 3: $0 < \al < 1$ and $\al + \be = 1$. Notice that in this case $k(\be-1)+\be \geq 0$ if $\alpha(k+1)<1$, so the
        best bound is obtained taking $k$ big enough such that $\al(k+1)>1$. In that case $\norm{T^n} = O(1)$.

        \item Case 4: $0 < \al < 1$, $0 < \be < 1$  and $\al + \be > 1$. Consider $\al = \frac{1}{m+1}$ with $m \in \NN$,
            then we have 3 possibilities. Take $k = m - 1$ with a growth of $O(n^{(\be-1)(m-1)+\be})$, take $k=m$ to achieve
            $O(n^{(\be-1)m+\be}\log(n))$ or take $k = m + 1$, and we have a growth of $O(n^{(\al + \be-1)(m+2)})$. Expressing $m$
            in terms of $\al$ we have that,
            \[
                n^{(\be-1)(m-1)+\be}    = n^{(\be-1)\left(\frac{1}{\al}-1\right)+\be} n^{1-\be}
                                            \geq  n^{(\be-1)\left(\frac{1}{\al}-1\right)+\be} \log(n)
                                            = n^{(\be-1)m+\be} \log(n),
            \]
            and that,
            \[
                n^{(\al+\be-1)(m+2)}   = n^{(\be-1)\left(\frac{1}{\al}-1\right)+\be} n^{\al+\be-1}
                                            \geq  n^{(\be-1)\left(\frac{1}{\al}-1\right)+\be} \log(n)
                                            = n^{(\be-1)m+\be} \log(n),
            \]
            for every $n$ sufficiently large. Thus, in this case the optimal integer to take is $k = \frac{1}{\al}-1$, and we obtain that
            $\norm{T^n} = O(n^{\frac{\al+\be-1}{\al}}\log(n))$. Consider now that $\alpha$ is not the inverse of a natural number. There
            are two options: take $k=\floor{\frac{1}{\al}-1}$ and the growth of the norm of the powers of $T$ is 
            $O(n^{(\be-1)\left(\floor{\frac{1}{\al}}-1\right)+\be})$, and the second option is to take $k = \floor{\frac{1}{\al}}$, and 
            we have that $\norm{T^n} = O(n^{(\al+\be-1)\left(\floor{\frac{1}{\al}}+1\right)})$. We introduce the following notation,
            we denote by $\eta = \frac{1}{\al} - \floor{\frac{1}{\al}}$, and now we study which $k$ is better,
            \[
                \begin{aligned}
                    (\al+\be-1)\left(\floor{\frac{1}{\al}}+1\right) - (\be-1)\left(\floor{\frac{1}{\al}}-1\right)-\be & =
                    \al \floor{\frac{1}{\al}} + \al + \be -2 \\ &= \al \left(\frac{1}{\al} - \eta \right) + \al + \be -2 \\
                    & = \al + \be - 1 -\al \eta.
                \end{aligned}                
            \]
            Therefore, if $\eta \leq \frac{\al+\be-1}{\al}$ then the optimal integer is $k=\floor{\frac{1}{\al}-1}$, and we have that
            $\norm{T^n} = O(n^{\floor{\frac{1}{\al}}(\be-1)+1})$. If $\eta \geq \frac{\al+\be-1}{\al}$ the optimal integer is
            $k = \floor{\frac{1}{\al}}$, and we have that $\norm{T^n} = O(n^{(\al+\be-1)\left(\floor{\frac{1}{\al}}+1\right)})$.
            In particular, we can see the representation of this different regions in \Cref{fig:regions_operator_powers}.

        \item Case 5: $0 < \al < 1$, $\be \geq 1$. First consider the case where $\al$ is not the inverse of a positive integer.
            Therefore, we just have to compare $O(n^\be)$ and $O(n^{(\al+\be-1)(k+1)})$ with $k = \floor{\frac{1}{\al}}$. Thus,
            \[
                (\al+\be-1)\left(\floor{\frac{1}{\al}}+1\right)-\be   > (\al+\be-1)\frac{1}{\al} - \be 
                                                                                = (\be - 1)\left(\frac{1}{\al}-1\right) > 0.
            \]
            So the best estimate is obtained for $k=0$, and it is $O(n^\be)$. On the other hand if $\al = \frac{1}{m+1}$ for 
            $m \in \NN$ we have to compare $O(n^\be)$, $O(n^{(\be-1)\left(\frac{1}{\al}-1\right) + \be} \log(n))$ and 
            $O(n^{(\al + \be-1)\left(\frac{1}{\al}+1\right)})$. Notice that, 
            \[
                n^{(\be-1)\left(\frac{1}{\al}-1\right)+\be} \log(n)  \geq  n^{\be} \log(n) \geq n^{\be},
            \]
            for $n$ sufficiently large. And lastly, as
            \[
                (\al+\be-1)\left(\frac{1}{\al}+1\right)-\be = \frac{1}{\al}(\be-1)+ \al > 0, 
            \]
            we have that 
            \[
                n^{(\al + \be-1)\left(\frac{1}{\al}+1\right)} \geq n^\beta.
            \]
            So, again the best estimate is obtained for $k=0$, and it is $O(n^\be)$.

        \item Case 6: $\al = 1$, $\be = 0$. Take $k=1$ and we have $\norm{T^n} = O(1)$.
        \item Case 7: $\al = 1$, $\be > 0$. Take $k=0$ and we have $\norm{T^n} = O(n^\be \log(n))$.
        \item Case 8: $\al > 1$, $\be \geq 1$. Take $k = 0$ and we have that $\norm{T^n} = O(n^{\al + \be -1})$.
    \end{itemize}
\end{proof}

Notice that in particular we recover the results seen in \Cref{re:exp_decay}. Once we have studied the growth of the norm of the powers,
our next aim is to study how the differences of these powers are bounded using similar techniques.

\begin{theorem}
    \label{th:differences_growth}
    Let $T$ be a $(\al,\be)$-RK operator on a complex Banach space $X$ with $\al,\be \geq 0$ and $k \in \NN \cup \{0\}$.
    Then, there exists $M_{\al,\be,k} \,>0$ depending on $\al,\be$ and $k$ such that,
    \begin{itemize}
        \item if  $\al(k+1) < 2$ then $\norm{T^{n+1}-T^n} \leq M_{\al,\be,k} \, n^{k(\be-1)+\be}$,
        \item if  $\al(k+1) = 2$ then $\norm{T^{n+1}-T^n} \leq M_{\al,\be,k} \, n^{k(\be-1)+\be} \log(n)$,
        \item if  $\al(k+1) > 2$ then $\norm{T^{n+1}-T^n} \leq M_{\al,\be,k} \, n^{(\al+\be-1)(k+1)-1}$.
    
    \end{itemize}
\end{theorem}

\begin{proof}
    We will prove the case for $k \in \NN$, for $k=0$ the steps are the same but is easier to compute. We start as in
    \Cref{th:powers_growth} using Cauchy differentiation formula of any order to the function  
    $f(\lambda) = \lambda^{n+k+1} - (1+\frac{k}{n+1}) \lambda^{n+k}$ to obtain that,
    \[
        T^{n+1}-T^n = \binom{n+k+1}{k}^{-1} \frac{1}{2\pi i} \int_{\Gamma} \lambda^{n+k}(\lambda-r) R(\lambda,T)^{k+1}\, d\lambda,
            \quad n \in \NN
    \]
    where $\Gamma = \{\lambda \in \CC: \abs{\lambda} = r\}$ with $r = 1+\frac{k}{n+1}>1$ if $k \in \NN$. For the case 
    $k=0$ we would have taken $r = 1+\frac{1}{n+1}>1$. Now, we introduce the bounds of the resolvent
    \begin{equation}
        \label{eq:bound_differences_powers_intermediate_formula}
        \begin{aligned}
            \norm{T^{n+1}-T^n}  & \leq \binom{n+k+1}{k}^{-1} \frac{C^{k+1} r^{n+ (\al+ \be)(k+1)-1}}{2\pi (r-1)^{\be(k+1)}}
                                    \int_{\Gamma} \frac{\abs{\lambda-r}}{\abs{\lambda-1}^{\al(k+1)}} d\lambda \\
                                & = \binom{n+k+1}{k}^{-1} \frac{C^{k+1} r^{n+ 1 + (\al+ \be)(k+1)}}{2\pi (r-1)^{\be(k+1)}}
                                    \int_{-\pi}^\pi \frac{\abs{e^{it}-1}}{\abs{re^{it}-1}^{\al(k+1)}} dt,
        \end{aligned}
    \end{equation}
    where the last integral is $J_{\al(k+1)}(r)$ studied in \Cref{lemma:integral_upper_bound}. The next step is to bound the terms
    appearing in the right-hand side of the last inequality. First, as in the proof of \Cref{th:powers_growth} we use the fact that
    $\binom{n+k+1}{k}^{-1} \leq  A_k (n+1)^{-k}  \leq  A'_k n^{-k}$ for some positive constant $A'_k$
    (see \cite[Equation 1.18, page 77, Vol.1]{Zygmund1968}). Secondly, as $r = 1+ \frac{k}{n+1}$ we have that
    $r^{n+1} = (1+\frac{k}{n+1})^{n+1} < e^k$. Also related with the value $r$ we have that
    $r^{(\al+\be)(k+1)} = (1+\frac{k}{n+1})^{(\al+\be)(k+1)}$  which is bounded by some constant as in the proof
    of \Cref{th:powers_growth}. Thus, \eqref{eq:bound_differences_powers_intermediate_formula} transforms into,
    \[
        \norm{T^{n+1}-T^n} \leq M_{\al,\be,k} n^{\be(k+1)-k} J_{\al(k+1)}(1+n^{-1}),
    \]
    where in $M_{\al,\be,k}$ are included all the constant that we have mentioned in the last paragraph. Lastly, we apply 
    the bounds seen in \Cref{lemma:integral_upper_bound}.
    \begin{itemize}
        \item If $\al(k+1) < 2$, then 
            \[
                \norm{T^{n+1}-T^n} \leq M_{\al,\be,k} \, n^{\be(k+1)-k} = M_{\al,\be,k} \, n^{(\be-1)k+\be}.
            \]
        \item If $\al(k+1) = 2$, then 
            \[
                \norm{T^{n+1}-T^n} \leq M_{\al,\be,k} \, n^{\be(k+1)-k} \log(n) = M_{\al,\be,k} \, n^{(\be-1) k+\be} \log(n).
            \]
        \item If $\al(k+1) > 2$, then 
            \[
                \norm{T^{n+1}-T^n} \leq M_{\al,\be,k} \, n^{\be(k+1)-k} n^{\al(k+1)-1} = M_{\al,\be,k} \, n^{(\al + \be-1)(k+1)}.
            \]
    \end{itemize}
    Where the constants of the \Cref{lemma:integral_upper_bound} enter into the constant $M_{\al,\be,k}$.
\end{proof}

Notice that with the case of Ritt and Kreiss operators we recover the known bounds of the power differences. For $\beta < 1$ and
$\alpha + \beta \geq 1$, we can achieve enhanced bounds for the norm differences of the operator $T$ powers. This improvement is based
on the results proved in \cite{Seifert2016} and further refinements in \cite{NgSeifert2020}. Note that in some sense these results are 
optimal.

\begin{theorem}
    \label{th:improved_differences_growth}
    Let $T$ be a $(\al,\be)$-RK operator on a Banach space $X$ with $\be<1$ and $\alpha + \beta \geq 1$, then 
    \[
        \norm{T^{n+1}-T^n} = O\left(\left( \frac{\log n}{n}\right)^{\frac{1-\be}{\al}}\right),
    \]
    for all sufficiently large $n$.
\end{theorem}
\begin{proof}
    By \Cref{pr:RK_resolvent_estimate} we have that 
    \[
        \norm{R(e^{i\theta},T)} \leq \frac{C_{\al,\be}}{\abs{e^{i\theta}-1}^\frac{\al}{1-\be}} 
                                = \frac{C_{\al,\be}}{2\sin \frac{\abs{\theta}}{2}}.
    \]
    Therefore, when $\theta \to 0$ we have that $R(e^{i\theta},T) = O(\abs{\theta})$, therefore by \cite[Corollary 3.1]{Seifert2016}
    the result follows.
\end{proof}

Observe that when $\alpha+\beta = 1$ in our \Cref{th:differences_growth} we can get rid of the logarithmic term to have that
$\norm{T^{n+1}-T^n} = O(n^{-1})$.

Now, we bring back the result proved in \cite[Theorem 4]{NagyZemanek1999} and commented in the introduction that characterizes the
Ritt operators through the growth of the following sequences $\{\norm{T^n}\}_{n \in \NN}$ and $\{ \norm{T^{n+1}-T^n}\}_{n \in \NN}$. 
In particular, we recall that a bounded linear operator $T$ is Ritt if and only if $\sup_{n \geq 0} \norm{T^n} < \infty$ and 
$\sup_{n \geq 0} n \norm{T^{n+1}-T^n} < \infty$. For our case we have the following corollary.

\begin{corollary}
    \label{cor:sum_implies_Ritt}
    Let $T$ be an $(\al,\be)$-RK operator with $\alpha+\beta=1$ and $\alpha>0$ then $T$ is Ritt.
\end{corollary}
\begin{proof}
    By \Cref{cor:regions_norms} and \Cref{th:differences_growth} $\norm{T^n} = O(1)$ and $\norm{T^{n+1} - T^n} = O(n^{-1})$ 
    and the result follows by the characterization of Ritt operators proven in \cite[Theorem 4]{NagyZemanek1999}.
\end{proof}

Note that we could have proved this result in differently way. As $\alpha+\beta=1$ and $\alpha>0$, from
\Cref{pr:RK_resolvent_estimate} we have that the Ritt's condition holds for $\TT \setminus \{1\}$. This already implies
by \cite{Lyubich1999} and \cite{LeMerdy2014} that the operator is Ritt. Although, we think that the proof presented using norm
estimates is more direct. As a consequence we have the following result.

\begin{corollary}
    There is no pair $(\al,\be)$ such that the set $\mathcal{RK}_{\al,\be}$ of operators satisfying the $(\al,\be)$-RK condition
    is the set of all power bounded operators.
\end{corollary}
\begin{proof}
    Denote by $PB$ the set of all power bounded operators and suppose $PB = \mathcal{RK}_{\al,\be}$. Since every power bounded
    operator is Kreiss and every Ritt is power bounded, 
    \[
        \mathcal{RK}_{1,0} \subset PB \subset \mathcal{RK}_{0,1},
    \]
    and by the inclusions seen in \eqref{eq:sets_RK} we have $\al+\be=1$. This yields, by \Cref{cor:sum_implies_Ritt},
    that if $\be<1$ every power bounded is Ritt, which is false, and if $\be=1$, then every Kreiss is power bounded, which is also
    false. Thus, the assertion is proved.
\end{proof}

\section{Applications, examples and counterexamples}
\label{sec:applications}

\subsection{ \texorpdfstring{$L^p$}{ Lᵖ} spaces and an interpolation result}

On of the main applications of this new condition is applying our results to operators on $L^p$ spaces with $1 \leq p \leq \infty$.
The study of growth of powers of operators on $L^p$ spaces has been widely studied, some references are 
\cite{DengLoristVeraar2023, Cuny2020,ArnoldCuny2023}.

Let $(\Omega, \mathcal{A}, \mu)$ be a measure space and $1 \leq p  \leq \infty$. Recall the spaces $L^p(\Omega)$ consists of 
equivalence classes of measurable functions $f: \Omega \rightarrow \mathbb{R}$ such that
\[
    \int \abs{f}^p \, d\mu < \infty, \; \text{for } \, 1 \leq p < \infty, \quad \text{ and } \quad 
    \esssup_{\Omega} \abs{f} < \infty \;\text{if } \, p = \infty,
\]
where two measurable functions are equivalent if they are equal $\mu$-almost everywhere. The $L^p$-norm of $f \in L^p(\Omega)$ 
is defined by
\[
    \|f\|_{p}=\left(\int|f|^p d \mu\right)^{1 / p}, \; \text{for } \, 1 \leq p < \infty, \quad \text{ and } \quad 
    \|f\|_{\infty}=  \esssup_{\Omega} \abs{f},      \; \text{if } \, p = \infty.
\]
Where the essential supremum of $f$ on $\Omega$ is,
\[
    \esssup_{\Omega} f  = \inf \{ a \in \RR: \mu \{ x \in \Omega : f(x) > a\}=0\}.
\]
From now on, we will write $L^p$ instead of $L^p(\Omega)$.

In particular, as a consequence of the Hölder's inequality applied to Ritt and Kreiss operators we have the following corollary.

\begin{corollary}
    \label{cor:interpolation_Ritt}
    Let $1 \leq p_0, p_1 \leq \infty$, and for $0 < \theta < 1$ define $p$  by
    \[
        \frac{1}{p} = \frac{1-\theta}{p_0}+\frac{\theta}{p_1}.
    \]
    If $T$ is a Kreiss operator from $L^{p_0}$ to $L^{p_0}$ and $T$ is a Ritt operator from $L^{p_1}$ to $L^{p_1}$, then
    $T$ is a Ritt operator from $L^p \to L^p$.
\end{corollary}
\begin{proof}
    As $T$ is a Kreiss operator from $L^{p_0}$ to $L^{p_0}$ we have that the resolvent operator $R(\lambda,T)$ is a bounded linear operator
    for $\abs{\lambda} > 1$ and
    \[
        \norm{R(\lambda , T)f}_{p_0} \leq \frac{C_1}{\abs{\lambda}-1} \norm{f}_{p_0}, \quad \text{for } f \in L^{p_0}.
    \]
    Analogously, as $T$ is a Ritt operator from $L^{p_1}$ to $L^{p_1}$ we have that,
    \[
        \norm{R(\lambda , T)f}_{p_1} \leq \frac{C_2}{\abs{\lambda-1}} \norm{f}_{p_1}, \quad \text{for } f \in L^{p_1}.
    \]
    Using Holder's  inequality  we have that $T$ is a bounded operator from $L^p \to L^p$ and that
    \[
        \norm{R(\lambda , T)f}_{p} \leq \frac{C_1^{1-\theta}C_2^\theta}{(\abs{\lambda}-1)^{1-\theta}(\abs{\lambda-1})^\theta} \norm{f}_{p},
        \quad \text{for } f \in L^{p} \text{ and } \abs{\lambda} > 1.
    \]
    This implies that $T$ is a $(\theta,1-\theta)$-RK operator. Thus, by \Cref{cor:sum_implies_Ritt} we have that $T$ is a Ritt operator
    from $L^p \to L^p$.
\end{proof}

\begin{remark}
    This result is important in the following sense. Usually when we study estimates about the growth of the norms in general Banach spaces
    this can be improved if we introduce the restriction to be working inside a Hilbert space. See for example \cite{BonillaMuller2021}
    for results on Hilbert spaces that improves the bounds known for Kreiss operators or \cite{MontesSánchezZemánek2015} where it is shown
    that the Volterra operator is power bounded if and only if we are on $L^2([0,1])$. Thus, \Cref{cor:interpolation_Ritt} brings 
    new ways to study this kind of operators as we can study the behavior of our operator on $L^2$ and then try to find what happens
     in other $L^p$ spaces.
\end{remark}

\subsection{ \texorpdfstring{$(C,\alpha)$}{(C, α)} bounded operators}
For $\alpha \in \CC$, we denote by $(k_\alpha(n))_{\NN_0}$ the sequence of Taylor coefficients of the generating function
$(1-z)^{-\alpha}$ where $\NN_0 := \NN \cup \{0\}$, that is,
\[
    \sum_{n=0}^\infty k_\alpha(n) z^n = \frac{1}{(1-z)^\alpha}, \quad \abs{z} < 1.  
\]
The elements of the sequences $(k_\alpha(n))_{\NN_0}$ are called Cesàro numbers, and are given by $(k_\alpha(0)) = 1$ and,
\[
    k_\alpha(n) = \frac{\Gamma(n+\alpha)}{\Gamma(\alpha) \Gamma(n+1)}, \quad \alpha \in \CC \setminus \{0,-1,-2,\dots\},
\]
where $\Gamma$ is the Gamma function. Moreover, let $T$ be a bounded linear operator and denote the discrete semigroup generated by 
the powers of the generator by $\mathcal{T}(n) := T^n$, $n \in \NN_0$, we define the Cesáro sum of order $\alpha \geq 0$ of 
$T$ as the family of bounded operators given by
\[
    S_n^{\al}\mathcal{T}(n)x:= \sum_{j=0}^{n} k^{\al}(n-j)T^{j}x,\, x\in X, n \in \mathbb{N}.
\]
And, the Cesáro mean of order $\al\geq 0$ of $T$ is the family of bounded operators given by
\[
    M_T^{\al}(n)x:= \frac{1}{k^{\al+1}(n)}S_n^{\al}\mathcal{T}(n)x,\, x\in X, N \in \NN_0.
\]
\begin{definition}
    We say $T$ is a $(C,\alpha)$ bounded operator if there exists $C>0$ such that 
    \[
        \sup_{n\geq 0}\|M_T^{\gamma}(n)\|\leq C.
    \]
\end{definition}

\begin{proposition}
    \label{pr:C_alpha_implies_RK}
    Let $T$ be a $(C,\alpha)$ bounded operator, then $T$ is an $(0,\alpha+1)$-RK operator.
\end{proposition}
\begin{proof}
    As $T$ is $(C,\alpha)$ bounded we have that there exists the algebra homomorphism $\theta$ and in that case by 
    \cite[Equation (2.6)]{AbadiasLizamaMianaVelasco2016} we have that
    \[
        \norm{R(\lambda,T)} = \frac{\abs{\lambda -1}^\alpha}{(\abs{\lambda} - 1)^{\alpha+1}}, \quad \abs{\lambda} > 1,
    \]
    thus, 
    \[
        \norm{R(\lambda,T)} \leq \frac{2^\al \abs{\lambda}^\alpha}{(\abs{\lambda} - 1)^{\alpha+1}}, \quad \abs{\lambda} > 1,
    \]
    which finishes the proof.
\end{proof}

Exploring Cesàro type operators and their power behavior is a common topic in literature
(refer to \cite{BermudezBonillaMullerPeris2020, BonillaMuller2021, CohenCunyEisnerLin2020} for more details).
However, a general issue arises when trying to identify $(C,\alpha)$ bounded operators as $(0,\alpha+1)$-RK operators.
The problem lies in the fact that the accuracy of this identification is limited, and it becomes suboptimal when
estimating the powers of these operators. For instance, if we consider a $(C,\alpha)$ bounded operator like $T$,
its norm $\norm{T^n}$ grows at a rate of $O(n^\alpha)$, as shown in \cite{Yoshimoto1998}. However, with our results,
we find that $\norm{T^n}$ exhibits a growth rate of $O(n^{\alpha+1})$, which is notably worse.

Another example from Bonilla and Muller in \cite[Example 1 and Example 4]{BonillaMuller2021} features a mean ergodic
operator $T$ on a Hilbert space $H$ that is not Kreiss but is clearly $(C,1)$. In this case, $\norm{T^n}$ grows at a
rate of $O(n)$; however, according to \Cref{pr:C_alpha_implies_RK} and \Cref{cor:regions_norms}, $T$ is classified
as an $(0,2)$-RK operator, resulting in a growth rate of $\norm{T^n} = O(n^2)$ a less favorable outcome than the
perspective of being $(C,1)$.

This isn't a new problem, as the traditional method of studying Cesàro type operators through their resolvent estimates
seems inherently problematic. Both Kreiss and Cesàro operators share the characteristic of $\norm{T^n} = O(n)$, but it's
clear that being Cesàro doesn't imply being Kreiss, as highlighted by the equality in the resolvent norm seen in the
proof of \Cref{pr:C_alpha_implies_RK}.

\section*{Acknowledgements}

We would like to extend our thanks to the referee for their valuable insights and constructive feedback, significantly improving
the quality of this paper. Also, we would like to thank L. Abadias, P.J. Miana, J. Oliva-Maza and A. Pritchard for the enriching
conversations and their helpful comments.

Version to appear in Studia Math.

\printbibliography 

\end{document}